\documentclass[12pt]{amsart}
\usepackage{amsmath,amssymb,amsbsy,amsfonts,latexsym,amsopn,amstext,
                                               amsxtra,euscript,amscd,bm,cite}

\usepackage{verbatim}

\usepackage{url}
\usepackage[colorlinks,linkcolor=blue,anchorcolor=blue,citecolor=blue,backref=page]{hyperref}

\usepackage{color}
\usepackage{float}
\usepackage{mathrsfs}

\usepackage{amsmath,amsthm,mathtools}

\usepackage{amssymb}

\hypersetup{breaklinks=true}

\usepackage[norefs,nocites]{refcheck}

\begin{document}

\newtheorem{theorem}{Theorem}
\newtheorem{lemma}[theorem]{Lemma}
\newtheorem{claim}[theorem]{Claim}
\newtheorem{cor}[theorem]{Corollary}
\newtheorem{prop}[theorem]{Proposition}
\newtheorem{definition}{Definition}
\newtheorem{question}[theorem]{Question}
\newtheorem{remark}[theorem]{Remark}
\newcommand{\hh}{{{\mathrm h}}}

\numberwithin{equation}{section}
\numberwithin{theorem}{section}
\numberwithin{table}{section}

\def\sssum{\mathop{\sum\!\sum\!\sum}}
\def\ssum{\mathop{\sum\ldots \sum}}
\def\iint{\mathop{\int\ldots \int}}

\def\squareforqed{\hbox{\rlap{$\sqcap$}$\sqcup$}}
\def\qed{\ifmmode\squareforqed\else{\unskip\nobreak\hfil
\penalty50\hskip1em\null\nobreak\hfil\squareforqed
\parfillskip=0pt\finalhyphendemerits=0\endgraf}\fi}

\newfont{\teneufm}{eufm10}
\newfont{\seveneufm}{eufm7}
\newfont{\fiveeufm}{eufm5}
%
%
\newfam\eufmfam
     \textfont\eufmfam=\teneufm
\scriptfont\eufmfam=\seveneufm
     \scriptscriptfont\eufmfam=\fiveeufm
%
%
\def\frak#1{{\fam\eufmfam\relax#1}}

\newcommand{\bflambda}{{\boldsymbol{\lambda}}}
\newcommand{\bfmu}{{\boldsymbol{\mu}}}
\newcommand{\bfxi}{{\boldsymbol{\xi}}}
\newcommand{\bfrho}{{\boldsymbol{\rho}}}

\newcommand{\bfalpha}{{\boldsymbol{\alpha}}}
\newcommand{\bfbeta}{{\boldsymbol{\beta}}}
\newcommand{\bfphi}{{\boldsymbol{\varphi}}}
\newcommand{\bfpsi}{{\boldsymbol{\psi}}}
\newcommand{\bftheta}{{\boldsymbol{\vartheta}}}

\def\fK{Frak K}
\def\fT{Frak{T}}

\def\fA{{Frak A}}
\def\fB{{Frak B}}
\def\fC{\mathfrak{C}}

\def \balpha{\bm{\alpha}}
\def \bbeta{\bm{\beta}}
\def \bgamma{\bm{\gamma}}
\def \blambda{\bm{\lambda}}
\def \bchi{\bm{\chi}}
\def \bphi{\bm{\varphi}}
\def \bpsi{\bm{\psi}}

\def\eqref#1{(\ref{#1})}

\def\vec#1{\mathbf{#1}}


\def\cA{{\mathcal A}}
\def\cB{{\mathcal B}}
\def\cC{{\mathcal C}}
\def\cD{{\mathcal D}}
\def\cE{{\mathcal E}}
\def\cF{{\mathcal F}}
\def\cG{{\mathcal G}}
\def\cH{{\mathcal H}}
\def\cI{{\mathcal I}}
\def\cJ{{\mathcal J}}
\def\cK{{\mathcal K}}
\def\cL{{\mathcal L}}
\def\cM{{\mathcal M}}
\def\cN{{\mathcal N}}
\def\cO{{\mathcal O}}
\def\cP{{\mathcal P}}
\def\cQ{{\mathcal Q}}
\def\cR{{\mathcal R}}
\def\cS{{\mathcal S}}
\def\cT{{\mathcal T}}
\def\cU{{\mathcal U}}
\def\cV{{\mathcal V}}
\def\cW{{\mathcal W}}
\def\cX{{\mathcal X}}
\def\cY{{\mathcal Y}}
\def\cZ{{\mathcal Z}}
\newcommand{\rmod}[1]{\: \text{mod} \: #1}

\def\cg{{\mathcal g}}

\def\vr{\mathbf r}

\def\e{{\mathbf{\,e}}}
\def\ep{{\mathbf{\,e}}_p}
\def\em{{\mathbf{\,e}}_m}

\def\Tr{{\mathrm{Tr}}}
\def\Nm{{\mathrm{Nm}\,}}

 \def\SS{{\mathbf{S}}}

\def\lcm{{\mathrm{lcm}}}
\def\ord{{\mathrm{ord}}}

\def\({\left(}
\def\){\right)}
\def\fl#1{\left\lfloor#1\right\rfloor}
\def\rf#1{\left\lceil#1\right\rceil}

\def\mand{\qquad \text{and} \qquad}

\newcommand{\commM}[1]{\marginpar{%
\begin{color}{red}
\vskip-\baselineskip 
\raggedright\footnotesize
\itshape\hrule \smallskip M: #1\par\smallskip\hrule\end{color}}}

\newcommand{\commI}[1]{\marginpar{%
\begin{color}{magenta}
\vskip-\baselineskip 
\raggedright\footnotesize
\itshape\hrule \smallskip I: #1\par\smallskip\hrule\end{color}}}

\newcommand{\commK}[1]{\marginpar{%
\begin{color}{blue}
\vskip-\baselineskip 
\raggedright\footnotesize
\itshape\hrule \smallskip K: #1\par\smallskip\hrule\end{color}}}




\hyphenation{re-pub-lished}

\mathsurround=1pt

\def\bfdefault{b}
\overfullrule=5pt

\def \F{{\mathbb F}}
\def \K{{\mathbb K}}
\def \N{{\mathbb N}}
\def \Z{{\mathbb Z}}
\def \Q{{\mathbb Q}}
\def \R{{\mathbb R}}
\def \C{{\mathbb C}}
\def\Fp{\F_p}
\def \fp{\mathfrak p}
\def \fq{\mathfrak q}

\def\ZK{\Z_K}

\def \xbar{\overline x}
\def\e{{\mathbf{\,e}}}
\def\ep{{\mathbf{\,e}}_p}
\def\eq{{\mathbf{\,e}}_q}

\title[primes represented by quartic polynomials on average]{On primes represented by quartic polynomials on average}

\date{\today}

\author[Kam Hung Yau]{Kam Hung Yau}

\address{Department of Pure Mathematics, University of New South Wales,
Sydney, NSW 2052, Australia}
\email{kamhung.yau@unsw.edu.au}

\begin{abstract}
We obtain an upper bound for the distribution of primes in the form $n^4+k$ up to $x$, averaged over $k$ with small square-full part. As a corollary, we show that for almost all $k$, there is an expected amount of  primes in the form $n^4+k$ up to $x$.
\end{abstract}


\maketitle

\section{Introduction}

Dirichlet proved that any linear polynomial with coprime coefficients take prime values infinitely often. The case of showing a polynomial of degree two or more  representing prime values infinitely often is still open,  although there are stunning partial results due to Kuhn~\cite{K} and later by Iwaniec~\cite{I}. Kuhn's result~\cite{K} can be extended to show that any irreducible polynomial $an^2 + bn + c$ with $a >0$ and $c$ is a product of at most three primes for infinitely many $n$. Iwaniec~\cite{I} improved Kuhn's result~\cite{K} for the particular polynomial $n^2+1$ by showing $n^2+1$ is a product of at most two primes for infinitely many $n$.  For polynomials of two variables, Friedlander \& Iwaniec~\cite{FI} proved a spectacular theorem stating that the polynomial $x^2 + y^4$  attain prime values infinitely often.  Later Heath-Brown~\cite{HB} and Heath-Brown \& Li~\cite{HL} showed the polynomials $x^3 + 2y^3$ and $a^2 + p^4$ ($p$ prime) takes prime values infinitely often respectively.

A conjecture of Bateman-Horn~\cite{BH} says that if $f$ is an irreducible polynomial in $\mathbb{Z}[x]$ satisfying $\gcd \{ f(n) : n \in \mathbb{Z} \}=1$ then
$$
\sum_{n \le x} \Lambda \left (f(n) \right ) \sim \prod_{p} \left (1 - \frac{n_p-1 }{p-1} \right ) x, 
$$
where $n_p$ is the number of solutions to the congruence 
\begin{align*}
f(n) \equiv 0 \bmod{p},  \mathrlap{\text{\qquad $n \in \mathbb{Z}/p\mathbb{Z}$.}}
\end{align*}
Here $\Lambda$ is the von-Mongoldt function.  See~\cite{BZ2} for a discussion on this topic. 
We note that in the special case of $f(n)=n^r+k$, we can write
$$
n_p = \sum_{j=0}^{\gcd(p-1,r)-1} \chi^j(-k),
$$
where $\chi$ is a multiplicative character of order $\gcd(p-1,r)$. The representation of $n_p$ for this specific $f$ is crucial for the following series of results.

Baier \& Zhao~\cite{BZ} showed that for fixed $A,B>0$, and $x^{r} (\log )^{-A} \le y \le x^r$, we have
$$
\sum_{\substack{ k \le y \\ \mu^2(k)=1 }}  \left | \sum_{n \le x} \Lambda(n^r +k) - \prod_{p>2} \left (1- \frac{ n_p-1  }{p-1} \right ) x \right |^2 \ll \frac{yx^2}{(\log x)^B},
$$
for $r=2$. Later, Foo \& Zhao~\cite{FZ} proved the case of $r=3$. We remark that in the speacial case of $r=2$, Baier \& Zhao~\cite{BZ3} widened  the region to $ x^{1/2 + \varepsilon}\le y \le x$ by instead using the dispersion method of Linnik~\cite{L}.

In this paper we are concerned with showing the case of $r=4$. As in all the previous approaches~\cite{BZ, FZ}, we use the circle method to estimate
$$
\sum_{n \le x} \Lambda(n^4+k) = \int_0^1 \sum_{m \le x^4 +k} \Lambda(m) e(\alpha m) \sum_{n \le x} e(-\alpha (n^4 +k)) \ \mathrm{d} \alpha.
$$
Next, we divide the interval $[0,1]$ into two sets, the major arc $\mathfrak{M}$ and the minor arc $\mathfrak{m}$, which corresponds to the main and error term respectively. The computation for the integral over the major arc $\mathfrak{M}$ and the corresponding error terms are done in sections~\ref{section: major arc}--\ref{section: error in major arc}. The computation for the integral over the minor arc $\mathfrak{m}$ are done in section~\ref{section: minor arc}.

\section{Results}

We write $f \ll g$ or  $f = O(g)$ to mean there exist a $C>0$ such that $f \le C g$. When $M \le m < 2M$ we write $m \sim M$. 

For any integer $n \ge 1$, we can write $n =\ell^2 m$ where $\mu^2(m)=1$. We denote the square-full part of $n$ to be $\kappa(n)=\ell^2$, and so  $n$ is square-free if and only if $\kappa(n)=1$.

\begin{theorem} \label{thm: average}
For any fixed $A, B , \varepsilon >0$, we have for $x^4(\log x)^{-A} \le y \le  x^4$, that
$$
\sum_{\substack{k \le y \\ \kappa(k) \le y^{1/2 - \varepsilon }}} \left | \sum_{n \le x} \Lambda(n^4+k) - \mathfrak{S}(k)x \right |^2 \ll \frac{yx^2}{(\log x)^B},
$$
where the singular series is given by
$$
\mathfrak{S}(k) = \prod_{p>2} \left (1 - \frac{n_p-1}{p-1} \right ),
$$
and $n_p$ is the number of solutions to $n^4+k \equiv 0 \bmod{p}$ in $\mathbb{Z}/ p\mathbb{Z}$.
\end{theorem}

Theorem~\ref{thm: average} implies the following result immediately.

\begin{cor} \label{cor}
For any fixed $A, B, C, \varepsilon >0$ and $\mathfrak{S}(k)$ as defined in Theorem~\normalfont{\ref{thm: average}}, we have for $x^4(\log x)^{-A} \le y \le x^4$ that
$$
\sum_{n \le x} \Lambda(n^4 +k) =  \mathfrak{S}(k) x  + O\left (\frac{x}{(\log x)^B} \right )
$$
for all $k$  up to $y$ such that $\kappa(k) \le y^{1/2 - \varepsilon}$ with at most $O \left  (y(\log x)^{-C} \right  )$ exceptions.
\end{cor}

\section{Preliminaries}

We recall a special case of the large sieve over quartic characters from~\cite[Theorem 1.2]{GZ}.

\begin{lemma} \label{lem: quartic large sieve}
Let $(a_m)_{m \in \mathbb{N}}$ be a sequence of complex numbers. Then
\begin{align*}
\sum_{ q \sim Q} \sideset{}{^\star} \sum_{\substack{ \chi  \bmod q \\ \chi^4 = \chi_0, \chi^2 \neq \chi_0 }} & \left |   \sum_{m \sim M} \mu^2(m) a_m \chi(m) \right  |^2 \\
& \ll (QM)^{\varepsilon}  ( Q^{5/4} +Q^{2/3}M ) \sum_{m \sim M} \mu^2(m) |a_m|^2,
\end{align*}
where the $\star$ on the sum over $\chi$ restricts the sum to primitive characters.
\end{lemma}

Next, we recall a bound for linear multiplicative character sum independently proved by P\'olya, and Vinogradov, see~\cite[Theorem 12.5]{IK}.

\begin{lemma}[P\'olya-Vinogradov] \label{lem: polya-vinogradov}
For any non-principle character $\chi$ modulo $q$, we have
$$
\left | \sum_{M < n \le M+N} \chi(n) \right | \ll q^{1/2} \log q.
$$
\end{lemma}

We state a  large sieve inequality  for number fields from~\cite[Theorem 1]{H}.

\begin{lemma} \label{lem: large sieve number field}
Let $K$ be a number field and $\mathfrak{r}$ denote an ideal in $K$. Suppose $u(\mathfrak{r})$ is a complex-valued function defined on the set of ideals in $K$. We have
$$
\sum_{\mathcal{N}(\mathfrak{f}) \le Q} \frac{\mathcal{N}(f)}{\Phi(\mathfrak{f})} \sideset{}{^*} \sum_{\chi  \bmod \mathfrak{f}} \left | \sum_{\mathcal{N}(\mathfrak{r}) \le z} u(\mathfrak{r}) \chi(\mathfrak{r}) \right |^2 \ll (z + Q^2) \sum_{\mathcal{N}(\mathfrak{r}) \le z} |u(\mathfrak{r})|^2,
$$
where $\mathcal{N}(\mathfrak{f})$ denotes the norm of the ideal $\mathfrak{f}$, $\Phi(\mathfrak{f})$ is Euler's totient function generalized to the setting of number fields, the $*$ over the summation over $\chi$ indicates that $\chi$ is a primitive character of narrow ideal class group modulo $\mathfrak{f}$ and the implicit constant depends on $K$.
\end{lemma}

We recall the Duality principle in~\cite[Theorem 228]{HLP}.

\begin{lemma}[Duality principle]  \label{lem: duality principle}
For a finite square $T=(t_{mn})$ matrix with entries in the complex numbers. The follow statements are equivalent:

For any complex sequence $(a_n)$, we have
$$
\sum_{m} \left | \sum_{n} a_n t_{mn} \right |^2 \le D \sum_{n} |a_n|^2.
$$

For any complex sequence $(b_n)$, we have
$$
\sum_{n} \left | \sum_{m   } b_m t_{mn} \right |^2 \le D \sum_{m} |b_m|^2.
$$
\end{lemma}

Recall the Perron formula from~\cite{D}.

\begin{lemma}[Perron] \label{lem: perron}
Suppose that $y \neq 1$ is a positive real number then for $c,T>0$ we have
$$
\frac{1}{2 \pi i} \int_{c-iT}^{c+iT} \frac{y^s}{s} \  \mathrm{d} s=
\begin{cases}
1 + O \left (y^{c} \min \{1,T^{-1} |\log y|^{-1} \} \right ) &\mbox{if $y>1$,} \\
O \left  (y^{c} \min \{1,T^{-1} |\log y|^{-1} \} \right )  & \mbox{otherwise.}
\end{cases}
$$
\end{lemma}

Recall the Weyl bound  from~\cite[Proposition 8.2]{IK}.

\begin{lemma}[Weyl] \label{lem: weyl}
If $f(x)= \alpha x^k+ \ldots + a_0$ is a polynomial with real coefficients  and $k \ge 1$ then
$$
\left | \sum_{n \le N}e(f(n)) \right | \le 2 N \left  \{N^{-k} \sum_{-N < \ell_1, \ldots, \ell_{k-1} <N} \min \left ( N , \frac{1}{\left \lVert \alpha k! \prod_{i=1}^{k-1} \ell_i \right \rVert} \right) \right \}^{2^{1-k}}.
$$
Here $\lVert x \rVert$ is the distance of $x$ to the nearest integer .
\end{lemma}

Recall a result from Mikawa~\cite{M}.

\begin{lemma} \label{lem: mikawa}
Let
$$
\mathfrak{J}(q,\Delta) = \sum_{\chi  \bmod q} \int_{N}^{2N} \left |  \sum_{t < n < t+q\Delta}^{\#} \Lambda(n) \chi(n) \right |^2 \mathrm{d} t
$$
where the $\#$ over the summation symbol means that if $\chi=\chi_0$, then $\chi(n) \Lambda(n)$ is replaced by $\Lambda(n)-1$. Let $\varepsilon,A,B>0$ be given. If $q \le (\log N)^B$ and $N^{1/5+\varepsilon} < \Delta < N^{1-\varepsilon}$, then we have
$$
\mathfrak{J}(q,\Delta) \ll_{\varepsilon,A,B} (q\Delta)^2 N(\log N)^{-A}.
$$
\end{lemma}

Recall  a result by Gallagher~\cite[Lemma 1]{G}.

\begin{lemma} \label{lem: gallagher}
Let $2 < \Delta < N/2$ and $N < N' < 2N$. For arbitrary complex sequence $(a_n)_{n \in \mathbb{N}}$, we have
$$
\int_{|\beta|< \Delta^{-1}} \left | \sum_{N < n < N'} a_n e(\beta n) \right |^2 d\beta \ll \Delta^{-2} \int_{N-\Delta/2}^N \left | \sum_{\max \{t,N \} < n < \min \{t+\Delta/2,N' \}} a_n \right |^2 \ \mathrm{d} t. 
$$
\end{lemma}

We recall a classical result from~\cite{H}.

\begin{lemma}[Bessel] \label{lem: bessel}
Let $\phi_1, \phi_2 , \ldots, \phi_R$ be orthonormal members of an inner product spaces $V$ over $\mathbb{C}$ and let $\xi \in V$. Then
$$
\sum_{r=1}^{R} |(\xi, \phi_r)|^2 \le (\xi, \xi).
$$
\end{lemma}

We recall a zero-free region for $L$-function from~\cite[Theorem 5.35]{IK}.

\begin{lemma} \label{lem: zero free region}
Let $K \backslash \mathcal{Q}$ be a number field, $\xi$ a Hecke Grossencharakter modulo $(\mathfrak{m}, \Omega)$ where $\mathfrak{m}$ is a non-zero integral ideal in $K$ and $\Omega$ is a set of real infinite places where $\xi$ is ramified. Let the conductor $\Delta = |d_K|N_{K/\mathcal{Q}} \mathfrak{m}$. There exists an absolute effective constant $c' >0$ such that the L-function $L(\xi, s)$ of degree $d=[K: \mathcal{Q}]$ has at most a simple real zero in the region
$$
\sigma > 1-\frac{c'}{d \log \Delta (|t| +3)}.
$$
The exceptional zero can occur only for a real character and it is strictly less than $1$.
\end{lemma}

\section{The major arc} \label{section: major arc}

The major arc is written as follows
$$
\mathfrak{M} =  \bigcup_{q \le Q_1} \bigcup_{\substack{a=1 \\ (a,q)=1}}^{q} I_{a,q},
$$
where
$$
I_{a,q} = \left [\frac{a}{q} - \frac{1}{qQ_2} , \frac{a}{q} + \frac{1}{qQ_2} \right ], \quad  Q_1 = (\log x)^{c_1},  \quad Q_2=x^{1-\varepsilon},
$$
for some suitable fixed $c_1 >0$. We remark that when $x$ is sufficiently large, $Q_2 > Q_1$ so that the intervals $I_{a,q}$ with $q \le Q_1$ are disjoint.

We write $\alpha \in \mathfrak{M}$ as
$$
\alpha = \frac{a}{q} + \beta \quad  \mbox{ with } \quad |\beta| \le \frac{1}{qQ_2}.
$$

Let 
$$
S_1(\alpha) = \sum_{m \le x^4+k} \Lambda(m)e(\alpha m) \quad \mbox{ and } \quad S_2(\alpha) = \sum_{n \le x} e  (-\alpha n^4  ).
$$
We will first estimate $S_1$. For convenience we set
\begin{equation} \label{eq: z}
z= x^4 +k.
\end{equation}

We have
\begin{align}
S_{1}(\alpha) & = \sum_{m \le z} \Lambda(m)e(am/q)e(\beta m) \nonumber \\
& = \sum_{\substack{m \le z \\ (m,q)=1}} \Lambda(m)e(a m/q)e(\beta m) +O \left (   \log^2 z  \right ). \label{eq: main term for S_1}
\end{align}

For $(am,q)=1$, we  recall the identity
\begin{equation} \label{eq: expression for e(am/q)}
e(am/q) = \frac{1}{\varphi(q)} \sum_{\chi \bmod q} \chi(am) \tau(\overline{\chi}),
\end{equation}
where
$$
\tau(\chi) = \sum_{n=1}^{q} \chi(n) e(an/q)
$$
is the well-known Gauss sum. Applying this identity, the main term of~\eqref{eq: main term for S_1} is transformed into
\begin{align*}
& \frac{1}{\varphi(q)} \sum_{\chi  \bmod q} \tau(\overline{\chi}) \chi(a) \sum_{m \le z} \chi(m) \Lambda(m) e(\beta m).
\end{align*}
Separating the main term corresponding to the principal character $\chi_0$ and using a property of Ramanujan sum, we obtain
\begin{align*}
& \frac{\mu(q)}{\varphi(q)} \left  (\sum_{m \le z} e(\beta m) + \sum_{m \le z} (\Lambda(m)-1)e(\beta m) \right  ) \\
&\quad  + \frac{1}{\varphi(q)} \sum_{\substack{\chi  \bmod q  \\ \chi \neq \chi_0 }} \tau(\overline{\chi}) \chi(a) \sum_{m \le z} \chi(m) \Lambda(m) e(\beta m) \\
& = \frac{\mu(q)}{\varphi(q)} \sum_{m \le z} e(\beta m) + \frac{1}{\varphi(q)} \sum_{\substack{\chi \bmod q   }} \tau(\overline{\chi}) \chi(a) \sum_{m \le z}^{\#} \chi(m) \Lambda(m) e(\beta m) \\
& = T_1(\alpha) + E_1(\alpha), \mbox{ say.}
\end{align*}
Here the $\#$ over the summation symbol means that if $\chi=\chi_0$, then $\chi(n) \Lambda(n)$ is replaced by $\Lambda(n)-1$.
Hence
$$
 S_1(\alpha)= T_1(\alpha) + E_1(\alpha) +O \left (\log^2 z \right ).
$$

Next we work on $S_2(\alpha)$. Recalling~\eqref{eq: expression for e(am/q)}, we get
\begin{align*}
S_2(\alpha) & = \sum_{ n \le x} e(-an^4/q)e(-\beta n^4 ) \\
&= \frac{1}{\varphi(q)} \sum_{\chi  \bmod q} \chi(-a) \tau(\bar{\chi})\sum_{n \le x} \chi^4(n) e(-\beta n^4)    \\
&=  \sum_{d|q} \frac{1}{\varphi(q_1^*)} \sum_{\chi   \bmod q_1^*} \chi(-ad^*) \tau(\overline{\chi}) \sum_{\substack{n \le x \\ (n,q)=1}} \chi^4(n^*) e(-\beta n^4),
\end{align*}
where
$$
q^* = \frac{q}{d}, \quad n^* = \frac{n}{d}, \quad d^* = \frac{d^3}{(d^3,q^*)}, \quad q_1^* = \frac{q^*}{(d^3,q^*)}.
$$
Therefore by separating the the main term corresponding to $\chi=\chi_0$, we obtain
\begin{align*}
S_2(\alpha) & = \sum_{d|q} \frac{1}{\varphi(q_1^*)} \sum_{ \substack{\chi  \bmod q_1^* \\ \chi^4 = \chi_0 }} \chi(-ad^*) \tau(\overline{\chi}) \sum_{\substack{n \le x \\ (n,q)=d}}  e(-\beta n^4) \\
& \quad + \sum_{d|q} \frac{1}{\varphi(q_1^*)} \sum_{ \substack{\chi  \bmod q_1^* \\ \chi^4 \neq \chi_0 }} \chi(-ad^*) \tau(\overline{\chi}) \sum_{\substack{n \le x \\ (n,q)=d }} \chi^4(n^*) e(-\beta n^4) \\
& = T_2(\alpha) + E_2(\alpha), \mbox{ say.}
\end{align*}

We consider the main term $T_2(\alpha)$. Observe that for $G= (\mathbb{Z} / q_1^* \mathbb{Z})^*$ and $G^4 = \{ g^4 : g \in G\}$, we have by orthogonality
\begin{align*}
\sum_{ \substack{\chi   \bmod q_1^* \\ \chi^4 = \chi_0 }} \chi(-ad^*) \tau(\overline{\chi}) & =  \sum_{\substack{b  \bmod q_1^* \\ (b,q_1^*)=1 }}  e(b/q_1^*) \sum_{\substack{\chi  \bmod q_1^* \\ \chi^4 = \chi_0}} \chi(-(ad^*)^3 b) \\
& =  \frac{ \# G}{  \# G^4} \sum_{\substack{-(ad^*)^3 b    \equiv \heartsuit  \bmod q_1^*  \\(b,q_1^*)=1}}  e(b/q_1^*),
\end{align*}
since $\chi^3 = \overline{\chi}$. Here $r \equiv \heartsuit  \bmod q_1^*$ means that $n$ is congruent to a fourth power of an integer modulo $q_1^*$. Moreover,
\begin{align*}
\sum_{\substack{-ad^*b \equiv \heartsuit  \bmod q_1^*  \\(b,q_1^*)=1}}  e(b/q_1^*) &  =  \frac{  \# G^4}{ \# G} \sum_{\substack{l=1 \\ (\ell,q_1^*)=1}}^{q_1^*} e(-\overline{(ad^*)^3} \ell^4/q_1^*) \\
& =\frac{  \# G^4}{ \# G} \sum_{\substack{\ell=1 \\ (\ell,q_1^*)=1}}^{q_1^*} e(-ad^* \ell^4/q_1^*),
\end{align*}
since  
$$
\overline{(ad^*)^3} \ell^4 \equiv ad^*  ( \overline{ad^{*} } \ell )^4     \bmod{q_1^*},
$$
where $\overline{u}$ is the multiplicative inverse of $u$ modulo $q_1^*$. It follows
$$
T_2(\alpha) = \sum_{d|q} \frac{1}{\varphi(q_1^*)} \sum_{\substack{\ell=1 \\ (\ell,q_1^*)=1}}^{q_1^*} e(-ad^* \ell^4/q_1^*)  \sum_{\substack{n \le x \\ (n,q)=d}}  e(-\beta n^4).
$$

In total, we obtain
\begin{align*}
&\int_{\mathfrak{M}} \sum_{m \le z} \Lambda(m) e(\alpha m)  \sum_{n \le x} e(-\alpha (n^4 +k)) \ \mathrm{d} \alpha \\
& = \int_{\mathfrak{M}} \{ T_1(\alpha) + E_1(\alpha) + O(\log^2 x) \} \{ T_2(\alpha) + E_2(\alpha) \} e(-\alpha k) \  \mathrm{d} \alpha.
\end{align*}

\section{Computing the singular series} 

The main term contributing in the major arc is
\begin{align*}
& \int_{\mathfrak{M}} T_1(\alpha) T_2(\alpha) e(-\alpha k) \ \mathrm{d}  \alpha \\
& = \sum_{q \le Q_1 } \frac{\mu(q)}{\varphi(q)} \sum_{\substack{a   \bmod q \\ (a,q)=1}} e \left (\frac{-ak}{q} \right ) \sum_{d|q} \frac{1}{\varphi(q_1^*)} \sum_{\substack{\ell =1 \\ (\ell, q_1^*)=1}}^{q_1^*} e \left (\frac{-ad^* \ell^4 }{q_1^*} \right  ) \int_{|\beta| < \frac{1}{qQ_2}} \Theta_{d,q}(\beta) \ \mathrm{d} \beta,
\end{align*}
where
$$
 \Theta_{d,q}(\beta) = e(-\beta k)  \sum_{\substack{n \le x \\ (n,q)=d }} e \left (-\beta n^4 \right ) \sum_{m \le z} e \left (-\beta m \right ). 
$$

As in~\cite{FZ}, we get that
\begin{align*}
& \int_{\mathfrak{M}} T_1(\alpha) T_2(\alpha) e(-k\alpha) \ \mathrm{d} \alpha \\
&= \sum_{q \le Q_1 } \frac{\mu(q)}{\varphi(q)} \sum_{\substack{a \bmod q \\ (a,q)=1}} e \left (\frac{-ak}{q} \right ) \sum_{d|q} \frac{1}{\varphi(q_1^*)} \sum_{\substack{\ell =1 \\ (\ell, q_1^*)=1}}^{q_1^*} e \left (\frac{-ad^* \ell^4 }{q_1^*} \right  ) \\
& \quad \times  \left \{ \frac{\varphi(q/d)}{q}x + O \left ( \left ( \frac{qQ_2x}{d}  \right )^{\frac{1}{2}} \right )  \right \}.
\end{align*}

The presence of $\mu(q)$ allows us to assume that $q$ is square-free and as a consequence we have $d^* =d^3$, $q_1^* = q^*=q/d$. Thus we further simplify the above to
\begin{align*}
& x\sum_{q \le Q_1 } \frac{\mu(q)}{\varphi(q)q} \sum_{\substack{a  \bmod q \\ (a,q)=1}} e \left (\frac{-ak}{q} \right ) \sum_{d|q}  \sum_{\substack{\ell =1 \\ (\ell, q/d)=1}}^{q/d} e \left (\frac{-a (d\ell)^4 }{q} \right  )  + O \left ((xQ_2)^{1/2} (\log x)^{c_2}  \right )
\end{align*}
for some fixed constant $c_2>0$.

Denote
\begin{align*}
\Sigma(q)  & = \sum_{\substack{a  \bmod q \\ (a,q)=1}} e \left (\frac{-ak}{q} \right ) \sum_{d|q}  \sum_{\substack{\ell =1 \\ (\ell, q/d)=1}}^{q/d} e \left (\frac{-a (d\ell)^4 }{q} \right  ) \\
& =  \sum_{\substack{a \bmod q \\ (a,q)=1}} e \left (\frac{-ak}{q} \right ) \sum_{\substack{m =1 }}^{q} e \left (\frac{-a m^4 }{q} \right  ).
\end{align*}

First we suppose $q=p$, a prime. A property of Ramanujan sum gives
\begin{align*}
\Sigma(p) &= \sum_{m=1}^{p} \sum_{\substack{ a  \bmod p \\ (a,p)=1}} e \left (\frac{-a(k+m^4)}{p} \right ) \\
& = (p-1)n_{k,p}-(p-n_{k,p}) \\
& = p(n_{k,p} -1 ),
\end{align*}
where $n_{k,p}$ is the number of solutions to 
\begin{align*}
m^4 + k \equiv 0  \bmod{p},  \mathrlap{\text{\qquad $m \in \mathbb{Z}/p\mathbb{Z}$.}}
\end{align*}
 Now let $\chi_{1,p}$, $\chi_{2,p}$, $\chi_{3,p}$ be three quartic characters of modulus $p$ prime such that $p\equiv 1 \bmod{4}$ and $\chi_{1,p}^2,\chi_{2,p}^2, \chi_{3,p}^2 \neq \chi_0$. Then
$$
\Sigma(p) =
\begin{cases}
p \{ \chi_{1,p}(-k) +\chi_{2,p}(-k) + \chi_{3,p}(-k) \} & \mbox{if $p\equiv 1 \bmod{4}$,} \\
0 & \mbox{if $p=2$ or $p \equiv 3  \bmod{4}$.}
\end{cases}
$$
Moreover, $\Sigma(q)$ is a multiplicative function in $q$.  Indeed for square-free positive integers $q_1,q_2$ with all primes factors congruent to 1 modulo 4, and $(q_1,q_2)=1$, we have
$$
\Sigma(q_1)\Sigma(q_2)= \sum_{m_1=1}^{p} \sum_{\substack{ a_1  \bmod p \\ (a_1,p)=1}}   \sum_{m_2=1}^{p} \sum_{\substack{ a_2   \bmod p \\ (a_2,p)=1}} e(\Upsilon_k(a_1,r_1,q_1, a_2,r_2,q_2)),
$$
where
\begin{align*}
\Upsilon_k(a_1,r_1,q_1, a_2,r_2,q_2) & = -k\frac{a_1 q_2 + a_2q_1}{q_1 q_2} - \frac{a_1q_2(q_2r_1)^4 + a_2q_1(q_1r_2)^4}{q_1q_2} \\
& \equiv -k\frac{a_1 q_2 + a_2q_1}{q_1 q_2} - \frac{(a_1q_2+a_2q_1)(q_1r_2+q_2r_1)^4}{q_1q_2}  \bmod{1}.
\end{align*}
Consequently
$$
 \Sigma(q_1 q_2) = \Sigma(q_1) \Sigma(q_2) .
$$
Therefore for some fixed constant $c_2 >0$, we get
\begin{align}
 \int_{\mathfrak{M}} T_1(\alpha) T_2(\alpha) e(-k\alpha) \ \mathrm{d} \alpha &  = x \sum_{q \le Q_1 } \frac{\mu(q)}{\varphi(q)q} \Sigma(q) +   O \left ((xQ_2)^{1/2} (\log x)^{c_2} \right ) \nonumber \\
 & = x \sum_{q =1}^{\infty} \frac{\mu(q)}{\varphi(q)} \prod_{p|q}(n_{k,p}-1)- x\Psi(k)  +   O \left ((xQ_2)^{1/2} (\log x)^{c_2}  \right )  \nonumber \\
 & = \mathfrak{S}(k) x + O \left ( x |\Psi(k)| +(xQ_2)^{1/2} (\log x)^{c_2} \right ), \label{eq: marjor arc}
\end{align}
where
$$
\Psi(k)=\sum_{q >Q_1} \frac{\mu(q)}{\varphi(q)} \prod_{p|q}(n_{k,p}-1).
$$
Note that we can restrict $q$ to be square-free with prime factors that are congruent to $1$ modulo $4$.

\section{Bounding the second moment of $\Psi(k)$} \label{se tion: second moment of Psi}

First we partition the second moment into three pieces
$$
\sum_{\substack{k \le y\\ \kappa(k) \le y^{1/2 -\varepsilon'} }} |\Psi(k)|^2 \ll \Psi_1 + \Psi_2 + \Psi_3,
$$
where
\begin{align*}
\Psi_1 & = \sum_{k \le y} \left | \sum_{Q_1 < q \le U }  \frac{\mu(q)}{\varphi(q)} \prod_{p|q}\{ \chi_{1,p}(-k) +\chi_{2,p}(-k) + \chi_{3,p}(-k) \} \right|^2, \\
\Psi_2 &=  \sum_{\substack{k \le y\\ \kappa(k) \le y^{1/2 -\varepsilon'} }}  \left | \sum_{U < q \le 2^vU }  \frac{\mu(q)}{\varphi(q)} \prod_{p|q}\{ \chi_{1,p}(-k) +\chi_{2,p}(-k) + \chi_{3,p}(-k) \} \right|^2, \\
\Psi_3 &=  \sum_{\substack{k \le y }}  \left | \sum_{ 2^v U < q }  \frac{\mu(q)}{\varphi(q)} \prod_{p|q}\{ \chi_{1,p}(-k) +\chi_{2,p}(-k) + \chi_{3,p}(-k)  \}  \right|^2,
\end{align*}
and $U,v$ are parameters to be chosen later.

 We deal with $\Psi_1$ first. Expanding the square, we obtain 
\begin{align*}
\Psi_1 \le &   \hspace{2mm} y \sum_{Q_1 < q \le U } \frac{\mu^2(q)}{\varphi(q)^2} 9^{\omega(q)}  + \sum_{ \substack{ Q_1 < q_1, q_2 \le U \\ q_1 \neq q_2 }} \frac{\mu(q_1) \mu(q_2)}{\varphi(q_1) \varphi(q_2) } \sum_{k \le y} P_{q_1,q_2}(k),
\end{align*}
where
\begin{align*}
& P_{q_1,q_2}(k) \\
&   = \prod_{p|q_1} \{ \chi_{1,p}(-k) +\chi_{2,p}(-k) + \chi_{3,p}(-k) \} \prod_{p|q_2}\{ \chi_{1,p}(-k) +\chi_{2,p}(-k) + \chi_{3,p}(-k) \}.
\end{align*}
The first term is no more than
$$
y \sum_{Q_1 < q \le U} \frac{9^{\omega(q)} (\log \log  10q)^2}{q^2} \ll y Q_1^{\varepsilon} \sum_{Q_1 < q \le U} \frac{1 }{q^2} \ll \frac{y}{Q_1^{1-\varepsilon}},
$$
by using the well-known bounds 
\begin{equation} \label{eg: phi bound}
\frac{q}{\log \log 10q} \ll \varphi(q),
\end{equation}
and
$$
 \omega(q) \ll \frac{\log q}{\log \log q}.
$$ 
The second term is bounded by
\begin{align*}
\sum_{\substack{Q_1 < q_1,q_2 \le U \\ q_1 \neq q_2 }} \frac{(q_1q_2)^{1/2} \log (q_1q_2) 3^{\omega(q_1) + \omega(q_2) }}{\varphi(q_1) \varphi(q_2)} & \ll (\log U) \left (\sum_{Q_1 < q \le U}  \frac{q^{1/2} 3^{\omega(q)}}{\varphi(q)} \right )^2\\
&  \ll U^{1+\varepsilon}.
\end{align*}
Therefore
\begin{equation} \label{eq: Psi_1 bound}
\Psi_1 \ll \frac{y}{Q_1^{1-\varepsilon}}  + U^{1+\varepsilon}.
\end{equation}

Partitioning the summation inside and applying the Cauchy's inequality, we have
$$
\Psi_2 \ll v \sum_{r=1}^{\lfloor v+1 \rfloor} \sum_{\substack{k \le y \\ \kappa(k) \le y^{1/2 -\varepsilon'} }} \left |  \sum_{2^{r-1}U < q \le 2^rU }  \frac{\mu(q)}{\varphi(q)} \prod_{p|q}\{ \chi_{1,p}(-k) +\chi_{2,p}(-k) + \chi_{3,p}(-k) \}  \right |^2.
$$
Expanding the term in the square, we have
\begin{align*}
 & \sum_{2^{r-1}U < q \le 2^rU }  \frac{\mu(q)}{\varphi(q)} \prod_{p|q}\{ \chi_{1,p}(-k) +\chi_{2,p}(-k) + \chi_{3,p}(-k) \}    \\
 & = \sum_{2^{r-1}U < q \le 2^rU } \sum_{\substack{\chi   \bmod q \\ \chi^4=\chi_0, \chi^2 \neq \chi_0}} \frac{\mu(q)}{\varphi(q)} \chi(k).
\end{align*}

Note that by Lemma~\ref{lem: quartic large sieve}, we get by Cauchy's inequality 
\begin{align*}
 & \sum_{2^{r-1}U < q \le 2^rU } \sum_{\substack{\chi  \bmod q \\ \chi^4=\chi_0, \chi^2 \neq \chi_0}}  \left | \sum_{\substack{k \le y \\ \kappa(k) \le y^{1/2 - \varepsilon' } }} a_k \chi(k) \right |^2 \\
 &= \sum_{2^{r-1}U < q \le 2^rU } \sum_{\substack{\chi \bmod q \\ \chi^4=\chi_0, \chi^2 \neq \chi_0}} \left |  \sum_{\ell^2 \le y^{1/2 - \varepsilon' }} \sum_{\substack{m \le y /\ell^2}} \mu^2(m) a_{\ell^2 m} \chi(\ell^2) \chi(m) \right |^2  \\
 & \ll y^{1/4 - \varepsilon'/2 }  \sum_{2^{r-1}U < q \le 2^rU } \sum_{\substack{\chi  \bmod q \\ \chi^4=\chi_0, \chi^2 \neq \chi_0}} \sum_{\ell^2 \le y^{1/2 - \varepsilon' } }    \left |  \sum_{\substack{m \le y /\ell^2}} \mu^2(m) a_{\ell^2 m} \chi(m) \right |^2   \\
 & \ll y^{1/4 - \varepsilon'/2 }   (2^r U y)^{\varepsilon} ( (2^rU)^{5/4} +  (2^rU)^{2/3}y) \sum_{\substack{ k \le y \\ \kappa(k) \le y^{1/2 - \varepsilon' } }}  |a_k|^2.
\end{align*}
Hence by the Duality principle (Lemma~\ref{lem: duality principle}), we have
\begin{align*}
\sum_{\substack{k \le y \\ \kappa(k) \le y^{1/2 - \varepsilon' } }}   &  \left | \sum_{2^{r-1}U < q \le 2^rU }  \sum_{\substack{\chi  \bmod q \\ \chi^4=\chi_0, \chi^2 \neq \chi_0}}  \frac{\mu(q)}{\varphi(q)}\chi(k) \right |^2  \\
&  \ll y^{1/4 - \varepsilon'/2 }  (2^r U y)^{\varepsilon} ( (2^rU)^{5/4} +  (2^rU)^{2/3}y)  \sum_{2^{r-1} U <q \le 2^r U} \frac{1}{\varphi(q)^2} \\
&  \ll y^{1/4 - \varepsilon'/2 }   (2^r U y)^{\varepsilon} ((2^r U)^{1/4} + (2^r U)^{-1/3}y).
\end{align*}

Summing over $r=1, \ldots, R=\lfloor \log_{2} ( y^{3 + \varepsilon}/U) \rfloor$ (logarithm base 2), we obtain
$$
\sum_{r=1}^{R} \sum_{\substack{k \le y \\ \kappa(k) \le y^{1/2 - \varepsilon' } }}   \left | \sum_{2^{r-1}U < q \le 2^rU } \sum_{\substack{\chi  \bmod q \\ \chi^4=\chi_0, \chi^2 \neq \chi_0}}  \frac{\mu(q)}{\varphi(q)}\chi(k) \right |^2 \ll y^{1/4 - \varepsilon' } \left ( y^{3/4 + 6\varepsilon} + \frac{y^{1+\varepsilon}}{U^{1/3-\varepsilon}} \right ).
$$

Note that primitive quartic character of conductor $q$ coprime to $4$ can be realised as quartic residue symbols $\left (\frac{m}{n} \right )_4$ for some square-free $n \in \mathbb{Z}[i]$, $n \equiv 1 \bmod{4}$ and not divisible by any rational primes with norm $q$. See~\cite{IR} for a thorough background.

For large values of $r$'s with $R< r \le \lfloor v+1 \rfloor$, it is enough to bound
\begin{equation} \label{eq: large sieve}
 \sum_{\substack{n \in \mathbb{Z}[i] \\ \mathcal{N}(n) \le y^2 \\ n \equiv 1  \bmod{4} \\ \kappa(n) \le y^{1/2 - \varepsilon' } }} \left | \sum_{\substack{\pi \in \mathbb{Q}[i] \\ 2^{r-1} < \mathcal{N}(\pi)=q \le 2^r U}} a_{\pi} \chi_{\pi}(n)  \right |^2,
\end{equation}
where $a_{\pi} = \mu(q)/\varphi(q)$, $\chi_{\pi}(n) = \left (\frac{n}{\pi} \right )_4$. Note that by the quartic reciprocity, we assert
$$
\chi_{\pi}(n) = \chi_n(\pi) (-1)^{\frac{\mathcal{N}(n)-1}{4} \frac{\mathcal{N}(\pi) -1}{4}},
$$
therefore in view of~\cite[Lemma 3.2]{GZ}, applying Lemma~\ref{lem: large sieve number field} to \eqref{eq: large sieve} and recalling~\eqref{eg: phi bound}, we majorise by
\begin{align*}
& \ll (y^4 + 2^r U) \sum_{\substack{\mathcal{N}(\pi)=q \\ 2^{r-1}U < q \le 2^r U}} \frac{1}{\varphi(q)^2 } \\
& \ll (y^4 + 2^r U) \frac{\log \log(2^r U) \log( 2^r U)}{2^r U}.
\end{align*}
Hence
$$
\sum_{\substack{k \le y \\ \kappa(k) \le y^{1/2 - \varepsilon' }  }}   \left | \sum_{2^{r-1}U < q \le 2^rU } \sum_{\substack{\chi   \bmod q \\ \chi^4=\chi_0, \chi^2 \neq \chi_0}}  \frac{\mu(q)}{\varphi(q)}\chi(k) \right |^2 \ll (y^4 + 2^r U) \frac{\log \log(2^r U) \log( 2^r U)}{2^r U}.
$$
Summing over $r$ from $R$ to $\lfloor v+1 \rfloor$, we have
\begin{align*}
&\sum_{r=R}^{\lfloor v+1 \rfloor} \sum_{\substack{k \le y  \\ \kappa(k) \le y^{1/2 - \varepsilon' } }}   \left | \sum_{2^{r-1}U < q \le 2^rU } \sum_{\substack{\chi   \bmod q \\ \chi^4=\chi_0, \chi^2 \neq \chi_0}}  \frac{\mu(q)}{\varphi(q)}\chi(k) \right |^2 \\
& \ll y^4 (v + \log U)\log (v + \log U )\sum_{r=R}^{\lfloor v+1 \rfloor} \frac{1}{2^r U} +    \sum_{r=R}^{\lfloor v+1 \rfloor}  (r + \log U)^2       \\
& \ll y^{1-\varepsilon}(v+\log U) \log(v + \log U)+ v^3 + v^2 \log U + v(\log U)^2 .
\end{align*}
Therefore 
\begin{align} 
\Psi_2  \ll v \Bigg ( y^{1/4 - \varepsilon'/2 }  \left ( y^{3/4 + 6\varepsilon} + \frac{y^{1+\varepsilon}}{U^{1/3-\varepsilon}} \right ) & + y^{1-\varepsilon}(v+\log U) \log(v + \log U)  \label{eq: Psi_2 bound} \\
&  +  v^3 + v^2 \log U + v(\log U)^2 \Bigg ). \nonumber
\end{align}

Lastly, we bound $\Psi_3$. For primes $p \equiv 1  \bmod{4}$, $p$ splits in $\mathbb{Q}[i]$ as $p= \pi_{p,1}\pi_{p,2}$.
Consider the function
\begin{align*}
f(s,k) & = \prod_{p \equiv 1  \bmod{4}}  \left (1 - \frac{\chi_{1,p}(-k) + \chi_{2,p}(-k) + \chi_{3,p}(-k)}{(p-1)p^s} \right ) \\
& = \prod_{p \equiv 1  \bmod{4}} \left (1 -\frac{n_{k,p}-1 }{(p-1)p^s} \right )  = \sum_q \frac{\mu(q)}{\varphi(q)q^s} \prod_{p|q} (n_{k,p}-1) \\
& = \prod_{\substack{p \equiv 1   \bmod{4} \\ p= \pi_{p,1} \pi_{p,2} }} \left  (1 - \frac{ \left (\frac{k}{\pi_{p,1}} \right )_4+ \left (\frac{k}{\pi_{p,2}} \right )_4 }{(p-1)p^s} \right ).
\end{align*}

Clearly if $s=0$ then $f(s,k)= \mathfrak{S}(k)$, and also $f(s,k)$ has no poles with $\Re(s)>0$. If we take 
$$
b_q = \frac{\mu(q)}{\varphi(q)} \prod_{p|q} (n_{k,p}-1), 
$$
then the Dirichlet series associated with $b_q$ is  $f(s,k)$.

Let us next consider the following Hecke $L$-function
\begin{align*}
L\left (s+1, \left (\frac{k}{\cdot} \right )_4 \right ) & = \prod_{\pi} \left (1 - \frac{ \left (\frac{k}{\pi}   \right )_4}{N(s)^{s+1} } \right )^{-1} \\
& = \left (1 - \frac{\left (\frac{k}{1-i}  \right )_4}{2^{s+1}} \right )^{-1} \prod_{p \equiv 3  \bmod{4}} \left (1 - \frac{ \left (\frac{k}{p} \right )_4}{p^{s+1} } \right )^{-1} \\
&\quad \times  \prod_{\substack{p \equiv 1 \bmod{4} \\ p=\pi_{p,1} \pi_{p,2}}} \left (1 - \frac{ \left (\frac{k}{\pi_{p,1}} \right )_4}{p^{s+1} } \right )^{-1}\left (1 - \frac{ \left  (\frac{k}{\pi_{p,2}} \right )_4}{p^{s+1} } \right )^{-1}.
\end{align*}

It is natural to approximate $f$ by $L^{-1}$. Denote
$$
h(s,k) = L\left (s+1, \left (\frac{k}{\cdot} \right )_4 \right ) f(s,k),
$$
so that we can write
$$
f(s,k) =  L^{-1}\left (s+1, \left (\frac{k}{\cdot} \right )_4 \right )h(s,k).
$$
It can be shown that $h$ is absolutely bounded for all $\Re(s)>-1/2 + \varepsilon$ for any fixed $\varepsilon>0$. Applying the Perron formula (Lemma~\ref{lem: perron}), we have
\begin{align*}
\sum_{y_1 \le q \le y_2} b_q & = \frac{1}{2\pi i} \int_{C-iT}^{C+iT} L^{-1}\left (s+1, \left (\frac{k}{\cdot} \right )_4 \right )h(s,k) \frac{y_2^s - y_1^s}{s} \ \mathrm{d} s \\
& \quad  +O \left ( \sum_{j=1}^2 \sum_q |b_q| \left (\frac{y_j}{q} \right )^C \min \left  \{1 , T^{-1} \left |\log \left (\frac{y_j}{q}  \right ) \right |^{-1} \right \} \right ).
\end{align*}
for any $C,T>0$. By applying Lemma~\ref{lem: zero free region} as in~\cite{FZ}, we can bound
\begin{equation} \label{eq: Psi_3 bound}
\Psi_{3} \ll y \exp \left (-c_3 \sqrt{2^v U} \right ),
\end{equation}
for some fixed constant $c_3>0$.

Setting 
$$
U =y^{3/4}, \quad  v = \log_2\left (  \frac{\exp(y^{\varepsilon/3})}{U} \right ), \quad  \varepsilon = \frac{\varepsilon '}{28}
$$
 (logarithm base 2) and recalling~\eqref{eq: Psi_1 bound},~\eqref{eq: Psi_2 bound}, and~\eqref{eq: Psi_3 bound}, we get
\begin{equation} \label{eq: Psi square average bound}
\sum_{\substack{k \le y \\ \kappa(k) \le y^{1/2 - \varepsilon'}}} |\Psi(k)|^2 \ll \frac{y}{(\log x)^{c_4} },
\end{equation}
for some fixed constant $c_4 >0$.

\section{Error terms from the major arc} \label{section: error in major arc}

There are three terms we need to bound, namely
\begin{align*}
\mathcal{E}_1 &=\sum_{k \le y} \left | \int_{\mathfrak{M}} T_1(\alpha) E_2(\alpha) e(-\alpha k)  \ \mathrm{d} \alpha \right |^2, \\
\mathcal{E}_2 & =\sum_{k \le y} \left | \int_{\mathfrak{M}} T_2(\alpha) E_1(\alpha) e(-\alpha k) \ \mathrm{d} \alpha \right |^2, \\
\mathcal{E}_3& = \sum_{k \le y} \left | \int_{\mathfrak{M}} E_1(\alpha) E_2(\alpha) e(-\alpha k) \ \mathrm{d} \alpha \right |^2.
\end{align*}

Now
\begin{align*}
\mathcal{E}_1  \ll \sum_{q < Q_1} \sum_{\substack{a \bmod q \\ (a,q)=1}} \Bigg (  &  \int_{|\beta| <\delta} \left | T_1(a/q+\beta)E_{2}(a/q +\beta) \right|^2 \ \mathrm{d} \beta   \\
& + \int_{\delta <|\beta| < \frac{1}{qQ_2} } \left | T_1(a/q+\beta)E_{2}(a/q +\beta) \right|^2 \ \mathrm{d} \beta \Bigg ).
\end{align*}

By the summing the geometric series we can bound by
$$
T_1(a/q+\beta) \ll \min(z,\beta^{-1}),
$$
and therefore
\begin{equation} \label{eq: mathcal{E}_1}
\mathcal{E}_1 \ll z^2 \sum_{q < Q_1} \sum_{\substack{a   \bmod q \\ (a,q)=1}} \left ( \int_{|\beta| <\delta} \left | E_{2}(a/q +\beta) \right|^2 \ \mathrm{d} \beta + \frac{1}{\delta^{2}}  \int_{\mathfrak{M}} \left | E_{2}(a/q +\beta) \right|^2 \ \mathrm{d} \beta \right ).
\end{equation}
First consider
$$
 \sum_{q < Q_1} \sum_{\substack{a  \bmod q \\ (a,q)=1}} \int_{|\beta| < \frac{1}{qQ_2} } \Omega_q(\beta) \ \mathrm{d} \beta, 
$$
where 
$$
\Omega_q(\beta) = \left |  \sum_{d|q} \frac{1}{\varphi(q_1^*)} \sum_{ \substack{\chi \bmod q_1^* \\ \chi^4 \neq \chi_0 }} \chi(-ad^*) \tau(\overline{\chi}) \sum_{\substack{n \le x \\ (n,q)=d }} \chi^4(n^*) e(-\beta n^4) \right |^2.
$$
Partition the summation over $n$ into dyadic intervals and let $N< n \le N' \le 2N \le x$. Applying the Cauchy's inequality, it is enough to bound
$$
(\log x)^{c_5} \sum_{q \le Q_1} \sum_{d|q} \sum_{\substack{ \chi   \bmod q_1^* \\ \chi^4 \neq \chi_0}} \int_{|\beta| \le \frac{1}{qQ_2}  } \left | \sum_{ \substack{ n \sim N' \\ (n,q)=d}} \chi^4(n^*) e(-\beta n^4) \right |^2 \ \mathrm{d} \beta, 
$$
which is no more than
$$
\frac{(\log x)^{c_5}}{Q_2^2} \sum_{q \le Q_1} \frac{1}{q^2} \sum_{d|q} \sum_{\substack{ \chi  \bmod q_1^* \\ \chi^4 \neq \chi_0}} \int_{-x}^{2x} \left | \sum_{\substack{\max \{ t,N \} < n \le \min \{t+qQ_2/2,N'\} \\  (n,q)=d}} \chi^4(n^*)  \right |^2 \ \mathrm{d}t
$$
by Lemma~\ref{lem: gallagher}. The character sum is
\begin{align*}
 \sum_{\substack{\max \{ t,N \} < n \le \min \{t+qQ_2/2,N'\} \\  (n,q)=d}} \chi^4(n^*) & = \sum_{\substack{\max \{ t,N \}/d < n* \le \min \{t+qQ_2/2,N'\}/d \\  (n^*,q^*)=1}} \chi^4(n^*) \\
 & = \sum_{\substack{\max \{ t,N \}/d < n* \le \min \{t+qQ_2/2,N'\}/d }} \chi'(n^*)  \\
 & \ll \sqrt{q_1^*} \log q_1^*,
\end{align*}
where $\chi'$ is the non-trivial character modulo $q_1^*$ induced by the character $\chi^2$ modulo $q_1^*$. The last line follows by P\'olya-Vinogradov (Lemma~\ref{lem: polya-vinogradov}). The first term in the upper bound of~\eqref{eq: mathcal{E}_1} can be estimated similarly. Therefore
$$
\mathcal{E}_1 \ll z^2 x \delta^2 (\log x)^{c_5} + \frac{x(\log x)^{c_6}}{\delta^2 Q_2^2}
$$
for some fixed constants $c_5, c_6 >0$.

Now by Lemma~\ref{lem: bessel}, we get
\begin{align*}
\mathcal{E}_2  & \ll \int_{\mathfrak{M}} |T_2(\alpha)   E_1(\alpha)|^2 \ \mathrm{d} \alpha \\
&  \ll \sup_{\alpha \in \mathfrak{M}} |T_2(\alpha)|^2 \int_{\mathfrak{M}} |E_1(\alpha)|^2 \ \mathrm{d} \alpha \\
& \ll x^2 \int_{\mathfrak{M}} |E_1(\alpha)|^2 \ \mathrm{d} \alpha.
\end{align*}
As exactly in~\cite{FZ}, we get using Lemma~\ref{lem: mikawa} to get
\begin{align} \label{eq: E_1 bound}
\int_{\mathfrak{M}} |E_1(\alpha)|^2 \ \mathrm{d} \alpha & \ll \sum_{q < Q_1} \frac{q}{\varphi(q)} (qQ_2)^{-2} \mathfrak{J} (q,Q_2/2) + Q_1^3 Q_2(\log x)^2  \nonumber \\
& \ll \sum_{q < Q_1} \frac{q}{\varphi(q)} z (\log z)^{-A},
\end{align}
for any $A>0$. Thus
$$
\mathcal{E}_2 \ll \frac{x^2 z}{(\log x)^{c_7} }.
$$

Again by Lemma~\ref{lem: bessel} and \eqref{eq: E_1 bound} we obtain
\begin{align*}
\mathcal{E}_3  & \ll \int_{\mathfrak{M}} |E_1(\alpha) E_2(\alpha)|^2 \ \mathrm{d} \alpha \\
& \ll \sup_{\alpha \in \mathfrak{M}} |E_2(\alpha)|^2 \int_{\mathfrak{M}} |E_1(\alpha)|^2 \ \mathrm{d} \alpha \\
&\ll \frac{x^2 z}{(\log x)^{c_8}}.
\end{align*}

Therefore by taking $\delta = 1/\sqrt{z}$, and collecting the bounds, we obtain
\begin{align}
\mathcal{E}_1 + \mathcal{E}_2 + \mathcal{E}_3 & \ll z^2 x \delta^2 (\log x)^{c_5} + \frac{x(\log x)^{c_6}}{\delta^2 Q^2} + \frac{x^2 z}{(\log x)^{c_7}} \nonumber \\
& \ll zx(\log x)^{c_5} + \frac{x^2 z}{(\log x)^{c_7}}. \label{eq: bound for e1 e2 e3}
\end{align}

\section{The  minor arc} \label{section: minor arc}

Lastly we estimate the minor arc
\begin{align*}
\sum_{\substack{k \le y \\ \kappa(k) \le y^{1/2-\varepsilon'} }}  & \left | \int_{\mathfrak{m} }  S_1(\alpha) S_2(\alpha) e(-\alpha k) \ \mathrm{d} \alpha \right |^2 \\
& = \sum_{\substack{k \le y \\ \kappa(k) \le y^{1/2-\varepsilon'} }} \left | \int_{\mathfrak{m} }  \sum_{m \le z} \Lambda(m)e(\alpha m) \sum_{n \le x} e(-\alpha(n^4 +k)) \ \mathrm{d} \alpha \right |^2.
\end{align*}
Here the minor arc is given by
$$
\mathfrak{m} = [1/Q_2,1+ 1/Q_2] \backslash \mathfrak{M}.
$$
By Lemma~\ref{lem: bessel}, the above is majorised by
$$
 \sup_{\alpha \in \mathfrak{m}} |S_2(\alpha)|^2 \int_0^1 |S_1(\alpha)|^2  \ \mathrm{d} \alpha \ll ( z \log z ) \sup_{\alpha \in \mathfrak{m}} |S_2(\alpha)|^2 .
$$
Now by Weyl shift (Lemma~\ref{lem: weyl}), we get
$$
|S_{2}(\alpha)|^8 \ll x^4  \sum_{-x < \ell_1,\ell_2,\ell_3 <x} \min \left \{ x, \frac{1}{\lVert 24 \alpha \ell_1 \ell_2 \ell_3  \rVert} \right \}.
$$

By Dirichlet's theorem, we have an approximation
$$
\left  | \alpha - \frac{a}{q}  \right | \le\frac{1}{48 x^3 q},
$$
where $(a,q)=1$ and $1 \le q \le 48x^3$. Since $\alpha \in \mathfrak{m}$, we further have that $q > Q_1$. Therefore for $-x < \ell_1, \ell_2, \ell_3 <x$, we have
$$
\left |24  \ell_1 \ell_2 \ell_3 \left (\alpha -  \frac{a}{q} \right ) \right| \le \frac{1}{2q},
$$
which implies
$$
\frac{1}{\lVert 24  \ell_1 \ell_2 \ell_3 \alpha \rVert } \le \frac{2}{\lVert 24  \ell_1 \ell_2 \ell_3 a/q \rVert}.
$$
In a first step
\begin{align*}
 \sum_{-x < \ell_1,\ell_2,\ell_3 <x} &  \min \left \{ x, \frac{1}{\lVert 24 \alpha \ell_1 \ell_2 \ell_3  \rVert} \right \} \\
 &  \ll x \sum_{\substack{-x < \ell_1, \ell_2 ,\ell_3 <x \\ q|24\ell_1 \ell_2 \ell_3}} 1 + \sum_{\substack{-x < \ell_1, \ell_2 ,\ell_3 <x \\ q \nmid24\ell_1 \ell_2 \ell_3}}  \frac{2}{\lVert 24  \ell_1 \ell_2 \ell_3 a/q \rVert}.
\end{align*}

Now let $q'=(q,24)$ then
\begin{align*}
\sum_{\substack{-x < \ell_1, \ell_2 ,\ell_3 <x \\ q|24\ell_1 \ell_2 \ell_3}} 1 & \ll \sum_{q_1 q_2 |q'} \sum_{\substack{-x < \ell_1 < x \\ q_1 | \ell_1 }}  \sum_{\substack{-x < \ell_2 < x \\ q_2 |\ell_2}} \sum_{\substack{-x < \ell_3 <x \\ \frac{q'}{q_1 q_2}| \ell_3 }}1 \\
& \ll \sum_{q_1 q_2 |q'} \left ( \frac{x}{q_1}+1 \right ) \left ( \frac{x}{q_2}+1 \right ) \left ( \frac{x}{q'/(q_1 q_2)}+1 \right ) \\
& \ll \tau(q')^2 \frac{x^3}{q'} + \tau(q')^2 x \\
& \ll x^3(\log x)^{-c_{9}},
\end{align*}
for some fixed constant $c_9>0$, and also
$$
\sum_{\substack{-x < \ell_1, \ell_2 ,\ell_3 <x \\ q \nmid24\ell_1 \ell_2 \ell_3}}  \frac{2}{\lVert 24  \ell_1 \ell_2 \ell_3 a/q \rVert} \ll \sum_{\substack{1 \le \ell_1, \ell_2 ,\ell_3 <x }} \frac{q}{\ell_1 \ell_2 \ell_3} \ll x^{3+\varepsilon}.
$$
Therefore
$$
\sup_{\alpha \in \mathfrak{m}} |S_2(\alpha)|^2 \ll \frac{x^2}{(\log x)^{c_{9}} },
$$
and hence recalling~\eqref{eq: z}, we obtain
\begin{equation} \label{eq: minor arc bound}
\sum_{\substack{k \le y \\ \kappa(k) \le y^{1/2 - \varepsilon'} }} \left | \int_{\mathfrak{m} }  S_1(\alpha) S_2(\alpha) e(-\alpha k) \ \mathrm{d} \alpha \right |^2  \ll \frac{x^6}{(\log x)^{c_9}}.
\end{equation}

\section{Proof of Theorem~\ref{thm: average}}

By Cauchy's inequality and recalling~\eqref{eq: marjor arc},~\eqref{eq: Psi square average bound},~\eqref{eq: bound for e1 e2 e3}, and~\eqref{eq: minor arc bound}, we get 
\begin{align*}
\sum_{\substack{k \le y \\ \kappa(k) \le y^{1/2 - \varepsilon'} }} \left | \sum_{n \le x} \Lambda(n^4+k) - \mathfrak{S}(k)x \right |^2  & \ll  \frac{yx^2}{(\log x)^{c_4}} + zx(\log x)^{c_5} + \frac{x^2 z}{(\log x)^{c_7}} + \frac{x^6}{(\log x)^{c_9}} \\
& \ll \frac{yx^2}{(\log x)^{c_4}}
\end{align*}
since $x^4(\log x)^{-A} \le y \le  x^4$ and by taking $c_4$ to be sufficiently large.

\section*{Acknowledgement}

The author thanks L. Zhao for the problem and many helpful conversations and comments, and also I. E. Shparlinski for helpful comments.
This work is supported by an Australian Government Research Training Program (RTP) Scholarship.


\begin{thebibliography}{12}

\bibitem{BZ} S. Baier, L. Zhao, {\it Primes in quadratic progressions on average.\/} Math. Ann. {\bf 338} (4) (2007), 963--982. 



\bibitem{BZ2} S. Baier, L. Zhao, {\it On primes represented by quadratic polynomials.\/} Anatomy of integers, 159--166, CRM Proc. Lecture Notes {\bf 46}, Amer. Math. Soc., 2008.


\bibitem{BZ3} S. Baier, L. Zhao, {\it On primes in quadratric progressions.\/}  Int. J. Number Theory {\bf 5} (6) (2009), 1017-–1035. 

\bibitem{BH} P. T. Bateman, R. A. Horn, {\it A heuristic asymptotic formula concerning the distribution of prime numbers.\/} Math. Comp. \textbf{16}, 1962, 363--367.

\bibitem{D} H. Davenport, {\it Multiplicative number theory.\/} Third edition. Graduate Texts in Mathematics, 74. Springer-Verlag, New York, 2000. 




\bibitem{FZ} T. Foo, L. Zhao, {\it On primes represented by cubic polynomials.\/} Math. Z. {\bf 274} (1-2) (2013), 323--340. 

\bibitem{FI} J. B. Friedlander, H. Iwaniec, {\it The polynomial $X^2 + Y^4$ captures its primes.\/} Ann. of Math. (2) {\bf 148} (3) (1998), 945--1040. 



\bibitem{G} P. X. Gallagher, {\it A large sieve density estimate near $\sigma=1.$\/} Invent. Math. \textbf{11}, 1970, 329--339. 

\bibitem{GZ} P. Gao, L. Zhao, {\it Large sieve inequalities for quartic characters.\/} Q. J. Math. {\bf 63} (4) (2012), 891--917. 

\bibitem{H} P. R. Halmos, {\it Finite-Dimensional Vector Spaces.\/} D. Van Nostrand, New York (1958).

\bibitem{HLP} G. H. Hardy, J. E. Littlewood, G. P\'olya, {\it Inequalities.\/} Cambridge University Press, Cambridge (1964).


\bibitem{HB} D. R. Heath-Brown, {\it Primes represented by $x^3+2y^3.$\/} Acta Math. {\bf 186} (1) (2001), 1--84.

\bibitem{HL} D. R. Heath-Brown, X. Li, {\it Prime values of $a^2+p^4.$\/} Invent. Math. {\bf 208} (2) (2017), 441--499.

\bibitem{H} M. N. Huxley, {\it The large sieve inequality for algebraic number fields.  \normalfont{II}. Means of moments of Hecke zeta-functions.\/} Proc. London Math. Soc. (3) \textbf{21} (1970), 108--128. 



\bibitem{IR} K. Ireland, M. Rosen, {\it A classical introduction to modern number theory.\/} Second edition. Graduate Texts in Mathematics {\bf 84}. Springer-Verlag, New York, 1990. 

\bibitem{I} H. Iwaniec, {\it The half dimensional sieve.\/} Acta Arith. {\bf 29} (1) (1976), 65--95.


\bibitem{IK} H. Iwaniec, E. Kowalski, {\it Analytic number theory.\/} American Mathematical Society Colloquium Publications {\bf 53}. American Mathematical Society, Providence, RI, 2004.

\bibitem{K} P. Kuhn, {\it \"Uber die Primteiler eines Polynomsw.\/} Proceedings of the International Congress of Mathematicians. Vol. 2 (Amsterdam, 1954), Erven P. Noordhoff N. V., Groningen; North-Holland, Amsterdam, 1954, 35--37.

\bibitem{L} Ju. V. Linnik, {\it The dispersion method in binary additive problems.\/} American Mathematical Society, Providence, R.I.

\bibitem{M} H. Mikawa, {\it On prime twins.\/} Tsukuba J. Math. {\bf 15} (1) (1991), 19--29. 




\end{thebibliography}
\end{document}